\documentclass[11pt]{amsproc}
\usepackage{amssymb,amsmath,amsthm,anysize}
\usepackage{longtable}
\usepackage{booktabs}
\usepackage{multirow,bigstrut,bigdelim,color}
\usepackage[shortlabels]{enumitem}
\usepackage{multicol}

\marginsize{2cm}{2cm}{2cm}{2cm}

\newtheorem{Thm}{Theorem}[section]

\newtheorem{La}[Thm]{Lemma}

\newenvironment{Prf}{\noindent\textbf{Proof.}}{\hfill $\Box$ \medskip}

\newcommand{\GL}{\textrm{GL}}

\newcommand{\SU}{\textrm{SU}}

\newcommand{\PGU}{\textrm{PGU}}

\newcommand{\Sp}{\textrm{Sp}}

\newcommand{\Sz}{\textnormal{Sz}}
\newcommand{\G}{\textnormal{G}}
\newcommand{\FF}{\textnormal{F}}
\newcommand{\B}{\textnormal{B}}

\setlist[description]{leftmargin=1ex,font=\normalfont\bfseries\space,style=nextline}

\begin{document}

\title[Generalised polygons with a point-primitive Suzuki or Ree group]{Generalised polygons admitting a point-primitive almost simple group of Suzuki or Ree type}
\author{Luke Morgan and Tomasz Popiel}

\address{
Centre for the Mathematics of Symmetry and Computation\\
School of Mathematics and Statistics\\
The University of Western Australia\\
35 Stirling Highway, Crawley, W.A. 6009, Australia.\newline
Email: \texttt{\{luke.morgan, tomasz.popiel\}@uwa.edu.au}.
}

\thanks{The first and second authors acknowledge the support of the Australian Research Council Discovery Project grants DP120100446 and DP140100416, respectively. 
We thank John Bamberg for helpful discussions, particularly regarding cases~(i) and~(ii) in Section~\ref{sec:2F4}.}

\subjclass[2010]{primary 51E12; secondary 20B15, 05B25}

\keywords{generalised hexagon, generalised octagon, primitive permutation group}

\maketitle

\begin{abstract}
\noindent
Let $G$ be a collineation group of a thick finite generalised hexagon or generalised octagon $\Gamma$. 
If $G$ acts primitively on the points of $\Gamma$, then a recent result of Bamberg et al.~shows that $G$ must be an almost simple group of Lie type.
We show that, furthermore, the minimal normal subgroup $S$ of $G$ cannot be a Suzuki group or a Ree group of type $^2\G_2$, and that if $S$ is a Ree group of type $^2\FF_4$, then $\Gamma$ is (up to point--line duality) the classical Ree--Tits generalised octagon.
\end{abstract}

\section{Introduction} \label{sec:intro}

A {\em generalised $d$-gon} is a point--line incidence geometry $\Gamma$ whose bipartite incidence graph has diameter $d$ and girth $2d$. 
If each point of $\Gamma$ is incident with at least three lines, and each line is incident with at least three points, then $\Gamma$ is said to be {\em thick}. 
By the well-known Feit--Higman Theorem~\cite{FeitHigman}, thick finite generalised $d$-gons exist only for $d \in \{2,3,4,6,8\}$. 
In the present paper, we are concerned with the cases $d=6$ (generalised {\em hexagons}), and $d=8$ (generalised {\em octagons}).

A {\em collineation} (or {\em automorphism}) of $\Gamma$ is a permutation of the point set of $\Gamma$, together with a permutation of the line set, such that the incidence relation is preserved  (equivalently, an automorphism of the incidence graph of $\Gamma$ that preserves the parts).
The only known thick finite generalised hexagons and octagons arise as natural geometries for certain exceptional groups of Lie type: $\text{G}_2(q)$ and $^3\text{D}_4(q)$ are collineation groups of generalised hexagons, and $^2\FF_4(q)$ acts on a generalised octagon. 
In each case, the action of the collineation group is primitive on both the points and the lines of $\Gamma$, and transitive on the {\em flags} of $\Gamma$, namely the incident point--line pairs. 
Each action is also point-distance-transitive --- that is, transitive on each set of ordered pairs of points at a given distance from each other in the incidence graph --- and line-distance-transitive. 
Buekenhout and Van Maldeghem~\cite{BvMdistance} showed that point-distance-transitivity implies point-primitivity for a thick finite generalised hexagon or octagon, and proved that there exist no point-distance-transitive examples other than the known {\em classical} examples. 
The existence of other point-primitive or flag-transitive (thick finite) generalised hexagons or octagons remains an open question.

Schneider and Van Maldeghem \cite{SvM} showed that a group $G$ acting point-primitively, line-primitively, and flag-transitively on a thick finite generalised hexagon or octagon must be an almost simple group of Lie type. 
That is, $S \le G \le \operatorname{Aut}(S)$, with $S$ a finite simple group of Lie type. 
Bamberg et~al.~\cite{PointPrimitiveGH&O} then showed that point-primitivity alone is sufficient to imply the same conclusion. 
We continue this work here, treating the families of Lie type groups which are of fixed rank and  fixed characteristic.

\begin{Thm} \label{thm:new}
Let $G$ be a point-primitive collineation group of a thick finite generalised hexagon or generalised octagon $\Gamma$, with $S \le G \le \operatorname{Aut}(S)$ for some nonabelian finite simple group $S$. 
Then $S$ is not a Suzuki group or a Ree group of type $^2\G_2$. 
Moreover, if $S$ is a Ree group of type $^2\FF_4$, then, up to point--line duality, $\Gamma$ is isomorphic to the classical Ree--Tits generalised octagon.
\end{Thm}

Theorem~\ref{thm:new} is proved in three sections:
the  Suzuki groups are  considered in Section~\ref{sec:Sz};  the small and large Ree groups are dealt with in Sections~\ref{sec:2G2} and~\ref{sec:2F4}, respectively.

\section{Preliminaries} \label{sec:pre}

Let us first collect some basic facts and definitions. 
If a finite generalised hexagon or octagon $\Gamma$ is thick, then there exist constants $s,t \ge 2$ such that each point (line) of $\Gamma$ is incident with exactly $t+1$ lines ($s+1$ points), and $(s,t)$ is called the {\em order} of $\Gamma$. 
If $\mathcal P$ denotes the  point set of $\Gamma$, then \cite[p.~20]{HendrikBook}
\begin{equation} \label{pts}
|\mathcal P| = \begin{cases} (s+1)(s^2t^2+st+1) & \text{if } \Gamma \text{ is a generalised hexagon,}\\
(s+1)(st+1)(s^2t^2+1) & \text{if } \Gamma \text{ is a generalised  octagon.}
\end{cases}
\end{equation}
Moreover, the integers $st$ and $2st$ are squares in the respective cases where $\Gamma$ is a generalised hexagon or generalised octagon.

\begin{La} \label{2part}
Let $\mathcal{P}$ be the point set of a thick finite generalised hexagon or generalised octagon $\Gamma$.
\begin{itemize}
\item[\textnormal{(i)}] If $2^a$ divides $|\mathcal{P}|$, where $a\ge 1$, then $|\mathcal{P}| > 2^{3a}$.
\item[\textnormal{(ii)}] If $\Gamma$ is a generalised hexagon and $3^a$ divides $|\mathcal{P}|$, where $a\ge 1$, then $|\mathcal{P}| > 3^{3a-4}$. 
\item[\textnormal{(iii)}] If $\Gamma$ is a generalised octagon and $2^a3^b$ divides $|\mathcal{P}|$, where $a\ge 0$ and $b\ge 1$, then $|\mathcal{P}| > 2^a3^{2b}$.
\end{itemize}
\end{La}

\begin{Prf}
Let $(s,t)$ be the order of $\Gamma$.

(i) First suppose that $\Gamma$ is a generalised hexagon. 
Since $s^2t^2+st+1$ is odd, $2^a$ must divide $s+1$. 
In particular, $s+1\ge 2^a$, and hence $s\ge 2^{a-1}$. 
Therefore, $|\mathcal{P}| > (s+1)s^2t^2 \ge 2^a(2^{a-1})^{2}2^2 = 2^{3a}$. 
Now let $\Gamma$ be a generalised octagon. 
Since $2st$ is a square, $st$ must be even, so $(st+1)(s^2t^2+1)$ is odd, and hence $2^a$ must divide $s+1$. 
Therefore, $|\mathcal{P}| > (s+1)s^3t^3 \ge 2^a(2^{a-1})^32^3 = 2^{4a} > 2^{3a}$.

(ii) Since $s^2t^2+st+1$ is not divisible by $9$, $s+1$ must be divisible by $3^{a-1}$. 
In particular, $s+1\ge 3^{a-1}$, and hence $s\ge 3^{a-2}$. 
Therefore, $|\mathcal{P}| > (s+1)s^2t^2 \ge 3^{a-1}(3^{a-2})^22^2 > 3^{3a-4}$. 

(iii) Since $2st$ is a square, $st$ is even, so $s^2t^2+1$ is divisible by neither $2$ nor $3$. 
Hence, $2^a3^b$ divides $(s+1)(st+1)$. 
In particular, $(s+1)(st+1) \ge 2^a3^b$. 
Let us say that $s+1$ is divisible by $3^c$, and that $st+1$ is divisible by $3^d$, where $c+d=b$. 
If $c\ge 1$, then $s > 3^{c-1/2}$; and if $d\ge 1$, then $st > 3^{d-1/2}$. 
Also, $t-1 = (st+1) - (s+1)t$ is divisible by $3^{\min\{c,d\}}$, so $t > 3^{\min\{c,d\}}$. 
If $c \ge d$, then $c\ge 1$ and $t > 3^d$, and hence $|\mathcal{P}| > (s+1)(st+1)(st)^2 \ge 2^a3^b (3^{c-1/2}3^d)^2 = 2^a3^{b+2(c+d)-1} = 2^a3^{3b-1} \ge 2^a3^{2b}$. 
If $d > c$, then $d \ge (b+1)/2$, so $|\mathcal{P}| > (s+1)(st+1)(st)^2 \ge 2^a3^b (3^{d-1/2})^2 \ge 2^a3^b (3^{b/2})^2 = 2^a3^{2b}$. 
\end{Prf}

Recall that a permutation group $G \le \operatorname{Sym}(\Omega)$ acts {\em primitively} on the set $\Omega$ if it acts transitively and preserves no nontrivial partition of $\Omega$, and that this is equivalent to the stabiliser $G_\omega$ of a point $\omega \in \Omega$ being a maximal subgroup of $G$. 
A maximal subgroup $M$ of an almost simple group $G$ with minimal normal subgroup $S$ is said to be a {\em novelty} maximal subgroup if $S \cap M$ is not maximal in $S$. 
Our notation is mostly standard: we write $D_{n}$ for a dihedral group of order $n$; $C_n$ denotes a cyclic group of order $n$; $[n]$ denotes an unspecified group of order $n$; and, for $q$ a prime power, $E_q$ denotes an elementary abelian group of order $q$.
For information about the Suzuki and Ree simple groups  of Lie type, we refer the reader to \cite{WilsonBook}, and the other references mentioned below. 

\section{Proof of Theorem~\ref{thm:new}: $S$ a Suzuki group} \label{sec:Sz}

We now adopt the hypothesis of Theorem~\ref{thm:new}, assuming additionally that $S$ is isomorphic to $\Sz(q)={}^2\B_2(q)$, where $q=2^m$ with $m$ odd and at least $3$. 
Then
\[
|S| = q^2(q^2+1)(q-1) = q^2(q+\sqrt{2q}+1)(q-\sqrt{2q}+1)(q-1).
\]
The outer automorphism group of $S$ is cyclic of order $m$. 
If we let $\sigma$ denote a generator of this group, then we have $G = S : \langle \sigma^j \rangle$ for some divisor $j$ of $m$. 
Let $\mathcal{P}$ be the point set of $\Gamma$,  and let $x \in \mathcal{P}$. 
Observe first that the stabiliser $G_x$ cannot contain $S$: if it did, then $G_x$ would have the form $S:K$ for some maximal subgroup $K$ of the cyclic group $\langle \sigma^j \rangle$, and hence $|G:G_x|=|\mathcal{P}|$ would be a prime, which is seen to be impossible upon inspection of \eqref{pts}.
Now, as explained in \cite[Section~7.3]{ColvaBook}, $G$ has no novelty maximal subgroups.
Therefore, $S_x = G_x \cap S$ is a maximal subgroup of $S$, so $S$ itself acts primitively on $\mathcal{P}$, and hence to prove the theorem we may assume that $G=S$. 
The maximal subgroups of $S$ are \cite[Table 8.16]{ColvaBook}, up to conjugacy,
\begin{itemize}
\item[(i)] $E_q . E_q . C_{q-1}$, 
\item[(ii)] $D_{2(q-1)}$, 
\item[(ii)] $C_{q\pm\sqrt{2q}+1} : C_4$, 
\item[(iv)] $\Sz(q_0)$, where $q=q_0^r$ with $r$ prime and $q_0 > 2$. 
\end{itemize}

\subsection{Case~(i)} \label{sz parabolic}

Suppose that $S_x \cong E_q . E_q . C_{q-1}$. 
Suzuki~\cite{Suzuki} showed that $S$ is 2-transitive in this action. 
Since $S$ preserves the incidence relation on $\Gamma$, and therefore distance in the incidence graph of $\Gamma$, we have  that the diameter of the incidence graph is at most three, a contradiction.

\subsection{Cases~(ii)--(iv)}

For the remaining cases, we apply Lemma~\ref{2part}(i).
If $S_x \cong D_{2(q-1)}$, then
\[
|\mathcal{P}| = |S:S_x| = \tfrac{1}{2} q^2(q^2+1) = 2^{2m-1}(2^{2m}+1) < 2^{4m},
\]
contradicting Lemma~\ref{2part}(i) with $a=2m-1$, which says that $|\mathcal{P}| > 2^{6m-3}$. 

If $S_x \cong C_{q\pm\sqrt{2q}+1}: C_4$, then
\[
|\mathcal{P}| = |S:S_x| = \tfrac{1}{4}q^2(q\mp+\sqrt{2q}+1)(q-1) = 2^{2m-2}(2^m\mp 2^{(m+1)/2} + 1)(2^m - 1) < 2^{4m-1},
\]
contradicting Lemma~\ref{2part}(i) with $a=2m-2$, which says that $|\mathcal{P}| > 2^{6m-6}$.

Finally, suppose that $S_x \cong \Sz(q_0)$, where $q=q_0^r$ with $r$ prime and $q_0 > 2$. 
Writing $q_0=2^\ell$, we have
\[
|\mathcal{P}| = |S:S_x| = 2^{2\ell(r-1)} \frac{(2^{2\ell r} + 1)(2^{\ell r} - 1)}{(2^{2\ell} + 1)(2^{\ell} - 1)} < 2^{5\ell(r-1)+2}, 
\]
contradicting Lemma~\ref{2part}(i) with $a=2\ell(r-1)$, which says that $|\mathcal{P}| > 2^{6\ell(r-1)}$.

\section{Proof of Theorem~\ref{thm:new}: $S$ a Ree group of type $^2\G_2$} \label{sec:2G2}

We now adopt the hypothesis of Theorem~\ref{thm:new} and assume that $S \cong {}^2\G_2(q)$, where $q=3^m$ with $m$ odd and at least $3$. 
Then
\[
|S| = q^3(q^3+1)(q-1) = q^3(q+\sqrt{3q}+1)(q-\sqrt{3q}+1)(q^2-1).
\]
Let $\mathcal{P}$ be the point set of $\Gamma$, and let $x\in \mathcal{P}$. 
The outer automorphism group of $S$ is cyclic (of order $m$), so, as in Section~\ref{sec:Sz}, we first deduce that $G_x$ is a maximal subgroup of $G$ not containing $S$. 
The maximal subgroups of $G$ were determined by Kleidman~\cite[Theorem~C]{Kleidman}. In particular, $G$ has no novelty maximal subgroups, so it suffices to prove the theorem in the case where $G=S$.
The maximal subgroups of $S$ are, up to conjugacy, 
\begin{itemize}
\item[(i)] $E_q.E_q .E_q.C_{q-1}$, 
\item[(ii)] $C_2 \times \text{PSL}_2(q)$,
\item[(iii)] $(E_4 \times D_{(q+1)/2}) : C_3$, 
\item[(iv)] $C_{q\pm\sqrt{3q}+1} : C_6$, 
\item[(v)] ${}^2\G_2(q_0)$, where $q=q_0^r$ with $r$ prime. 
\end{itemize}

\subsection{Case~(i)}

Suppose that $S_x \cong E_q.E_q .E_q.C_{q-1}$. 
Then $S$ acts 2-transitively on $\mathcal{P}$ \cite[p.~251]{DixonMortimer}. 
The same argument as in Section~\ref{sz parabolic} now provides a contradiction.

\subsection{$\Gamma$ a generalised hexagon: cases~(ii)--(v)} \label{sec4.2}

For cases~(ii)--(v) with $\Gamma$ a generalised hexagon, we use Lemma~\ref{2part}(ii).
First suppose that $S_x \cong C_2\times \text{PSL}_2(q)$. 
The order of $S_x$ is $q(q^2-1)$, so
\[
|\mathcal{P}| = |S:S_x| = q^2(q^2-q+1) = 3^{2m}(3^{2m}-3^m+1) < 3^{4m},
\]
contradicting Lemma~\ref{2part}(ii) with $a=2m$, which says that $|\mathcal{P}| > 3^{6m-4}$. 

If $S_x \cong (E_4\times D_{(q+1)/2}) : C_3$, then
\[
|\mathcal{P}| = |S:S_x| = \tfrac{1}{6}q^3(q-1)(q^2-q+1) = \tfrac{1}{2} 3^{3m-1}(3^m-1)(3^{2m}-3^m+1) <  3^{6m-1},
\]
contradicting Lemma~\ref{2part}(ii) with $a=3m-1$, which says that $|\mathcal{P}| > 3^{9m-7}$. 

If $S_x \cong C_{q\pm \sqrt{3q} + 1} : C_6$, then
\[
|\mathcal{P}| = |S:S_x| = q^3(q^2-1)(q\mp \sqrt{3q} + 1) = 3^{3m}(3^{2m}-1)(3^m\mp 3^{(m+1)/2} + 1) < 3^{6m+1}, 
\]
contradicting Lemma~\ref{2part}(ii) with $a=3m$, which says that $|\mathcal{P}| > 3^{9m-4}$. 

Finally, suppose that $S_x \cong {}^2G_2(q_0)$, where $q=q_0^r$ with $r$ prime. 
Writing $q_0=3^\ell$, we have
\[
|\mathcal{P}| = |S:S_x| = 3^{3\ell(r-1)} \frac{(3^{3\ell r} + 1)(3^{\ell r} - 1)}{(3^{3\ell} + 1)(3^{\ell} - 1)} < 3^{7\ell(r-1)+2}.
\]
If $\ell(r-1) \ge 3$, then this contradicts Lemma~\ref{2part}(ii) with $a=3\ell(r-1)$, which gives $|\mathcal{P}| > 3^{9\ell(r-1)-4}$. 
Otherwise, $(\ell,r)=(1,3)$, and there is no valid solution $(s,t)$ to equation~\eqref{pts}.

\subsection{$\Gamma$ a generalised octagon: cases~(ii)--(iv)}

Now suppose that $\Gamma$ is a generalised octagon. 
We first use Lemma~\ref{2part}(iii) to rule out cases~(ii)--(iv) for $S_x$, computing $|S:S_x|$ in each case as in Section~\ref{sec4.2}.
First suppose that $S_x \cong C_2\times \text{PSL}_2(q)$. 
Then
\[
|\mathcal{P}| = |S:S_x| = 3^{2m}(3^{2m}-3^m+1) < 3^{4m}, 
\]
contradicting Lemma~\ref{2part}(iii) with $a=0$ and $b=2m$, which says that $|\mathcal{P}| > 3^{4m}$.
 
Next, suppose that $S_x \cong (E_4\times D_{(q+1)/2}) : C_3$. 
Observe that $3^{3m}+1$ is divisible by $4$, because $3m$ is odd. 
Therefore,
\[
|\mathcal{P}| = |S:S_x| = 2\cdot 3^{3m-1} \frac{3^{3m}+1}{4} < 2\cdot 3^{6m-2},
\]
while Lemma~\ref{2part}(iii) with $a=1$ and $b=3m-1$ gives $|\mathcal{P}| > 2\cdot 3^{6m-2}$, a contradiction.

Finally, suppose that $S_x \cong C_{q\pm \sqrt{3q} + 1} : C_6$. 
Observe that $3^{2m}-1$ is divisible by $2^3$ because $m$ is odd, and that $3^m\mp 3^{(m+1)/2} + 1$ is even. 
Therefore, 
\begin{align*}
|\mathcal{P}| = |S:S_x| &= 2^4 3^{3m} \frac{(3^{2m}-1)(3^m\mp 3^{(m+1)/2} + 1)}{2^4} \\
&\le 2^4 3^{3m} \frac{(3^{2m}-1)(3^m + 3^{(m+1)/2} + 1)}{2^4} < 2^4 3^{6m-2},
\end{align*}
while Lemma~\ref{2part}(iii) with $a=4$ and $b=3m$ gives $|\mathcal{P}| > 2^43^{6m}$, a contradiction.

\subsection{$\Gamma$ a generalised octagon: case~(v)}

Finally, we consider case~(v) with $\Gamma$ a generalised octagon. 
The approach is similar to that used for cases~(ii)--(iv), but requires a little more care.

Suppose that $S_x \cong {}^2\G_2(q_0)$, where $q=q_0^r$ with $r$ prime. 
Writing $q_0=3^\ell$, we have
\begin{equation} \label{smallReeBound}
|\mathcal{P}| = 3^{3\ell(r-1)} \frac{(3^{3\ell r} + 1)(3^{\ell r} - 1)}{(3^{3\ell} + 1)(3^{\ell} - 1)} < 3^{7\ell(r-1)+\epsilon}, 
\quad \text{where } \epsilon := \frac{\log\left(\frac{3^4}{(3^3-1)(3-1)}\right)}{\log(3)} \approx 0.336.
\end{equation}
To verify the inequality in \eqref{smallReeBound}, one checks that $(3^{3\ell} + 1)(3^\ell -1) \ge 3^{4\ell-\epsilon}$, because $\ell \ge 1$, and that $(3^{3\ell r} + 1)(3^{\ell r} - 1) < 3^{4\ell r}$. 
Let us re-write this inequality as
\[
|\mathcal{P}| < 3^{ 7b/3+\epsilon}, \quad \text{where } b := 3\ell(r-1).
\]
Note also that $b\ge 6$, because $r\ge 3$.
For a contradiction, we now show that $|\mathcal{P}| > 3^{7b/3+\epsilon}$. 
By \eqref{smallReeBound}, $3^b$ is the highest power of $3$ dividing $|\mathcal{P}|$. 
Since $2st$ is a square, $st$ is even, so $s^2t^2+1$ is not divisible by $3$. 
Hence, by \eqref{pts}, $3^b$ divides $(s+1)(st+1)$. 
As in the proof of Lemma~\ref{2part}(iii), let us say that $s+1$ is divisible by $3^c$, and that $st+1$ is divisible by $3^d$, where $b=c+d$. 
Recall also (from that proof) that $t > 3^{\min\{c,d\}}$. 
To show that $|\mathcal{P}| > 3^{7b/3+\epsilon}$, we consider four cases.

First suppose that $c \ge d$. 
Then $t > 3^d$, and $c \ge 1$ so $s \ge 3^c-1 > 3^{c-1/2}$. 
Hence, $|\mathcal{P}| > (s+1)(st+1)(st)^2 > 3^b (3^{c-1/2}3^d)^2 = 3^{b+2(c+d)-1} = 3^{3b-1} > 3^{7b/3+1}$, with the final inequality holding because $b \ge 6 > 3$. 
Next, suppose that $d/2+1/2 \le c < d$. 
Then $6 \le b = c+d \le 3c-1$. 
In particular, $c \ge (b+1)/3$; and $c\ge 3$ so $s \ge 3^c-1 \ge 3^{c-\delta}$, where $\delta := 3 - \log(3^3-1)/\log(3)$.
Moreover, $t>3^c$, and hence $|\mathcal{P}| > 3^b(st)^2 > 3^b(3^{c-\delta}3^c)^2 = 3^{b+4c-2\delta} \ge 3^{7b/3 + (4/3-2\delta)}$. 
It follows that $|\mathcal{P}| > 3^{7b/3 + \epsilon}$, because $1.26 \approx 4/3-2\delta > \epsilon \approx 0.336$.
Now suppose that $c \le d/2-1/2$. 
Then $6 \le b = c+d \le 3d/2-1/2$. 
In particular, $d\ge (2b+1)/3$; and $d\ge 5$ so $st \ge 3^d-1 \ge 3^{d-\delta'}$, where $\delta' := 5 - \log(3^5-1)/\log(3)$. 
Therefore, $|\mathcal{P}| > 3^b(st)^2 > 3^{b+2d-2\delta'} = 3^{7b/3+2/3-2\delta'}$, and it follows that $|\mathcal{P}| > 3^{7b/3 + \epsilon}$, because $0.659 \approx 2/3-2\delta' > \epsilon \approx 0.336$.

Finally, suppose that $d/2-1/2 < c < d/2+1/2$. 
Since $c$ and $d$ are integers, this is equivalent to saying that $c=d/2$. 
Now, suppose first, towards a contradiction, that $(s+1)(st+1)$ is actually equal to $3^b$. 
Then $s+1=3^c$, $st+1=3^{2c}$, and \eqref{smallReeBound} implies that
\begin{equation} \label{c=d/2}
(s^2t^2+1)(3^{3\ell} + 1)(3^{\ell} - 1) = (3^{3\ell r} + 1)(3^{\ell r} - 1).
\end{equation}
However, this is impossible, because the left- and right-hand sides of \eqref{c=d/2} are not congruent modulo $3$. 
Indeed, $st = 3^{2c}-1 \equiv 2 \pmod 3$, so $s^2t^2+1 \equiv 4+1 \equiv 2 \pmod 3$; $3^{3\ell} + 1 \equiv 1 \pmod 3$; and $3^\ell - 1 \equiv 2 \pmod 3$; and hence the left-hand side of \eqref{c=d/2} is congruent to $1$ modulo $3$. 
On the other hand, the right-hand side of \eqref{c=d/2} is congruent to $2$ modulo $3$.
Therefore, $(s+1)(st+1)$ is strictly larger than $3^b$. 
Indeed, it is larger by a factor of at least $5$, because by \eqref{smallReeBound} we see that $|\mathcal{P}|/3^b$ is divisible by neither $2$ nor $3$ (to verify that $|\mathcal{P}|/3^b$ is odd, apply \cite[Lemma~2.5]{GuestPraeger}). 
Therefore, $|\mathcal{P}| > 5\cdot 3^b(st)^2 > 3^{b+1}(st)^2$. 
Since $6 \le b = 3d/2$, we have $d\ge 4$, and so $st \ge 3^d-1 \ge 3^{d-\delta''}$, where $\delta'' := 4 - \log(3^4-1)/\log(3)$. 
Hence, $|\mathcal{P}| > 3^{b+1+2d-2\delta''} = 3^{7b/3+1-2\delta''}$, and it follows that $|\mathcal{P}| > 3^{7b/3 + \epsilon}$, because $0.977 \approx 1-2\delta''  > \epsilon \approx 0.336$.

\section{Proof of Theorem~\ref{thm:new}: $S$ a Ree group of type $^2\FF_4$} \label{sec:2F4}

In this final section, we adopt the hypothesis of Theorem~\ref{thm:new} while assuming that $S \cong {}^2\FF_4(q)$, where $q=2^m$ with $m$ odd and at least $3$. 
Then
\[
|S| = q^{12}(q^6+1)(q^4-1)(q^3+1)(q-1).
\]
Let $\mathcal{P}$ be the point set of $\Gamma$, and let $x\in \mathcal{P}$. 
The outer automorphism group of $S$ is cyclic, so we again observe that $G_x$ is a maximal subgroup of $G$ not containing $S$. 
A result of Malle~\cite{Malle} tells us that $G$ has no novelty maximal subgroups, so it again suffices to prove the theorem in the case where $G=S$. 
The maximal subgroups of $S$ (listed also in \cite[Section~4.9.3]{WilsonBook}) are, up to conjugacy,
\begin{multicols}{2}
\begin{itemize}
\item[(i)] $P_1 := [q^{10}] : (\Sz(q) \times C_{q-1})$,
\item[(ii)] $P_2 := [q^{11}] : \mathrm{GL}_2(q)$;
\item[(iii)] $\SU_3(q) : C_2$, 
\item[(iv)] $\PGU_3(q) : C_2$,
\item[(v)] $\Sz(q) \wr C_2$,
\item[(vi)] $\Sp_4(q) : C_2$,
\item[(vii)] ${}^2\FF_4(q_0)$, where $q = q_0^r$ with $r$ prime, 
\item[(viii)] $(C_{q + 1} \times C_{q+1}) : \GL_2(3)$,
\item[(ix)] $C_{(q \pm \sqrt{2q} + 1)^2} : [96]$,
\item[(x)] $C_{q^2 +q+1\pm \sqrt{2q}(q+1)} : C_{12}$.
\end{itemize}
\end{multicols}
The groups $P_1$ and $P_2$ are maximal parabolic subgroups of $S$. 
The group $P_1$ is a point stabiliser in the action of $S$ on the classical generalised octagon, whilst $P_2$ is a point stabiliser in the action of $S$ on the dual \cite[Section~4.9.4]{WilsonBook}. 
We must show that $S_x$ cannot be isomorphic to any of the groups in cases~(iii)--(x), and, further, that if $S_x$ is isomorphic to either $P_1$ or $P_2$, then $\Gamma$ is the classical generalised octagon or its dual.

\subsection{Cases~(i)--(ii) with $\Gamma$ a generalised octagon} \label{2F41stCase}

Suppose that $\Gamma$ is a generalised octagon and that $S_x$ is isomorphic to either $P_1$ or $P_2$. 
In either action, the group $S$ has rank five. 
That is, the point stabiliser $S_x$ has five orbits on the set $\mathcal P$ \cite[p.~167]{WilsonBook}. 
For $i\in\{0,2,4,6,8\}$, denote by $\Gamma_i(x)$ the set of points at distance $i$ from $x$ in the incidence graph of $\Gamma$. 
Since each of these sets is nontrivial and $S_x$-invariant, the pigeonhole principle shows that each is an orbit of $S_x$. 
Since $S$ acts transitively on $\mathcal P$, we find that $S$ acts distance-transitively on $\mathcal{P}$.
Now the main result of \cite{BvMdistance} shows that $\Gamma$ is isomorphic to the classical generalised octagon associated with $S$, or its dual.

\subsection{Case~(i) with $\Gamma$ a generalised hexagon} \label{large ree hexagon}

Suppose that $\Gamma$ is a generalised hexagon, with $S_x \cong [q^{10}] : (\Sz(q) \times C_{q-1})$. 
Since $|\Sz(q)| = q^2(q^2+1)(q-1)$,
\[
|\mathcal{P}| = |S:S_x| = (q^4-q^2+1)(q^3+1)(q^2+1)(q+1).
\]
Equivalently (subtracting $1$ from both sides),
\begin{equation} \label{2F4P1pts-1}
s^3t^2 + s^2(t+1) + s(t+1) = q^{10} + q^9 + q^7 + q^6 + q^4 + q^3 + q.
\end{equation}
Now, $S$ acts primitively and distance-transitively on the points of a generalised octagon of order $(q,q^2)$, with point stabiliser $[q^{10}] : (\Sz(q) \times C_{q-1})$ and nontrivial subdegrees \cite[Section~4.9.4]{WilsonBook}
\begin{equation} \label{subdegrees2F4points}
n_1 := q(q^2+1),\quad n_2 := q^4(q^2+1),\quad n_3 := q^7(q^2+1),\quad n_4 := q^{10}.
\end{equation}
Recall the notation $\Gamma_i(x)$ from Section~\ref{2F41stCase}. 
Then we have \cite[p.~19]{HendrikBook}
\begin{equation} \label{distancesHexagon}
|\Gamma_2(x)| = s(t+1),\quad |\Gamma_4(x)| = s^2t(t+1),\quad |\Gamma_6(x)| = s^3t^2,
\end{equation}
and $S_x$  preserves the sets $\Gamma_i(x)$.
Hence, each $\Gamma_i(x)$ is a union of $S_x$-orbits, and so for $i\in\{2,4,6\}$, we have $|\Gamma_i(x)| = \sum_{k=1}^4 \delta_{i,k} n_k$, for some $\delta_{i,k} \in \{0,1\}$ (with $\delta_{i,k}\delta_{j,k}=0$ for $i \neq j$). 
We show that this leads to a contradiction. 

{\bf Claim~$1$:} $|\Gamma_2(x)| = n_1$. 
The proof of the claim is by contradiction. 
If not, then $|\Gamma_2(x)| \ge n_2 = q^4(q^2+1)$. 
Since $s,t \ge 2$, and so in particular $t \ge \tfrac{2}{3}(t+1)$, it follows that
\begin{align*}
|\Gamma_4(x)| &\ge \tfrac{2}{3}s^2(t+1)^2 = \tfrac{2}{3} |\Gamma_2(x)|^2 \ge \tfrac{2}{3} q^8(q^2+1)^2, \\
|\Gamma_6(x)| &\ge 2s^2t^2 \ge \tfrac{4}{3} |\Gamma_4(x)| \ge \tfrac{8}{9} q^8(q^2+1)^2.
\end{align*}
Since the left-hand side of \eqref{2F4P1pts-1} is $|\Gamma_2(x)| + |\Gamma_4(x)| + |\Gamma_6(x)|$, this implies that
\[
\tfrac{14}{9} q^8(q^2+1)^2 + q^4(q^2+1) \le q^{10} + q^9 + q^7 + q^6 + q^4 + q^3 + q,
\]
which is certainly false for $q \ge 8$.

{\bf Claim~$2$:} $|\Gamma_4(x)| = n_2$. 
The proof is again by contradiction. 
If not, then $|\Gamma_4(x)| \ge n_3 = q^7(q^2+1)$, because $|\Gamma_2(x)| = n_1 = q(q^2+1)$ by Claim~1. 
This implies the following inequality, which is certainly false for $q \ge 8$:
\[
q^6 = \frac{q^7(q^2+1)}{q(q^2+1)} \le \frac{|\Gamma_4(x)|}{|\Gamma_2(x)|} = \frac{s^2t(t+1)}{s(t+1)} = st < s(t+1) = q(q^2+1).
\]

By Claims~1 and~2, we must have $|\Gamma_6(x)| = n_3+n_4 = q^7(q^3+q^2+1) > q^8(q^2+1)$, and hence 
\[
s > \frac{s^3t^2}{s^2t(t+1)} = \frac{|\Gamma_6(x)|}{|\Gamma_4(x)|} > \frac{q^8(q^2+1)}{q^4(q^2+1)} = q^4.
\]
This is impossible, because $s(t+1) = q(q^2+1)$ by Claim~1 (and hence certainly $s < q(q^2+1) < q^4$).

\subsection{Case~(ii) with $\Gamma$ a generalised octagon}

Suppose that $\Gamma$ is a generalised hexagon, with $S_x \cong [q^{11}] : \GL_2(q)$. 
Since $|\GL_2(q)| = q(q^2-1)(q-1)$,
\[
|\mathcal{P}| = |S:S_x| = (q^4-q^2+1)(q^2+1)^2(q^3+1).
\]
Equivalently (subtracting $1$ from both sides),
\begin{equation} \label{2F4P1pts-1ii}
s^3t^2 + s^2(t+1) + s(t+1) = q^{11} + q^9 + q^8 + q^6 + q^5 + q^3 + q^2.
\end{equation}
Now, $S$ acts primitively and distance-transitively with stabiliser $[q^{11}] : \GL_2(q)$ on the points of a generalised octagon of order $(q^2,q)$, namely the point--line dual of the generalised octagon from case~(i). 
The nontrivial subdegrees are \cite[p.~167]{WilsonBook}
\begin{equation} \label{subdegrees2F4pointsii}
n_1 := q^2(q+1),\quad n_2 := q^5(q+1),\quad n_3 := q^8(q+1),\quad n_4 := q^{11}.
\end{equation}
For $x\in \mathcal{P}$, we again have \eqref{distancesHexagon}, and $S_x$ must preserve the sets $\Gamma_i(x)$, $i\in\{2,4,6\}$, so each $|\Gamma_i(x)|$ is equal to a sum of   the subdegrees $n_1,\dots,n_4$, as in Section~\ref{large ree hexagon}.
We show that this leads to a contradiction. 

{\bf Claim~$1$:} $|\Gamma_2(x)| = n_1$. 
The proof of the claim is by contradiction. 
If not, then $|\Gamma_2(x)| \ge n_2 = q^5(q+1)$. 
Since $s,t \ge 2$, and so in particular $t \ge \tfrac{2}{3}(t+1)$, it follows that
\begin{align*}
|\Gamma_4(x)| &\ge \tfrac{2}{3}s^2(t+1)^2 = \tfrac{2}{3} |\Gamma_2(x)|^2 \ge \tfrac{2}{3} q^{10}(q+1)^2, \\
|\Gamma_6(x)| &\ge 2s^2t^2 \ge \tfrac{4}{3} |\Gamma_4(x)| \ge \tfrac{8}{9} q^{10}(q+1)^2.
\end{align*}
Since the left-hand side of \eqref{2F4P1pts-1ii} is $|\Gamma_2(x)| + |\Gamma_4(x)| + |\Gamma_6(x)|$, this implies the following inequality, which is false for $q \ge 8$:
\[
\tfrac{14}{9} q^{10}(q+1)^2 + q^5(q+1) \le q^{11} + q^9 + q^8 + q^6 + q^5 + q^3 + q^2,
\]

{\bf Claim~$2$:} $|\Gamma_4(x)| = n_2$. 
The proof is again by contradiction. 
If not, then $|\Gamma_4(x)| \ge n_3 = q^8(q+1)$, because $|\Gamma_2(x)| = n_1 = q^2(q+1)$ by Claim~1. 
This implies the following inequality, which is false for $q \ge 8$:
\[
q^6 = \frac{q^8(q+1)}{q^2(q+1)} \le \frac{|\Gamma_4(x)|}{|\Gamma_2(x)|} = \frac{s^2t(t+1)}{s(t+1)} = st < s(t+1) = q^2(q+1).
\]

By Claims~1 and~2, we must have $|\Gamma_6(x)| = n_3+n_4 = q^8(q^3+q+1) > q^9(q^2+1)$, and hence 
\[
s > \frac{s^3t^2}{s^2t(t+1)} = \frac{|\Gamma_6(x)|}{|\Gamma_4(x)|} > \frac{q^9(q^2+1)}{q^5(q+1)} = \frac{q^4(q^2+1)}{q+1}.
\]
This, however, contradicts $s(t+1) = q^2(q^2+1)$ (namely Claim~1).

\subsection{Cases (iii)--(ix)}

We now  deal with cases~(iii)--(ix), for which we use Lemma~\ref{2part}(i) to contradict the  equality $|\mathcal{P}| = |S:S_x|$.
First suppose that $S_x$ is isomorphic to either $\SU_3(q) : C_2$ or $\PGU_3(q) : C_2$. 
In either case, we have $|S_x| = 2q^3(q^3+1)(q^2-1)$, and hence 
\[
|\mathcal{P}| = |S:S_x| = \tfrac{1}{2} q^9(q^6+1)(q^2+1)(q-1) = 2^{9m-1}(2^{6m}+1)(2^{2m}+1)(2^m-1) < 2^{18m+1}.
\]
However, Lemma~\ref{2part}(i) with $a=9m-1$ gives $|\mathcal P | > 2^{27m-3}$, which is a contradiction.

If $S_x \cong \Sz(q) \wr C_2$, then $|S_x| = 2q^4(q^2+1)^2(q-1)^2$, so
\[
|\mathcal{P}| = |S:S_x| = \tfrac{1}{2}q^8(q^4-q^2+1)(q^3+1)(q+1) = 2^{8m-1}(2^{4m}-2^{2m}+1)(2^{3m}+1)(2^m+1) < 2^{16m+1},
\]
contradicting Lemma~\ref{2part}(i) with $a=8m-1$, which gives $|\mathcal{P}| > 2^{24m-3}$.
 
If $S_x \cong \mathrm{Sp}_4(q) : C_2$, then $|S_x| = 2q^4 (q^4-1)(q^2-1)$, so
\[
|\mathcal{P}| = |S:S_x| = \tfrac{1}{2} q^8(q^6+1)(q^2-q+1) = 2^{8m-1}(2^{6m}+1)(2^{2m}-2^m+1) < 2^{16m},
\]
contradicting Lemma~\ref{2part}(i) with $a=8m-1$, which gives $|\mathcal{P}| > 2^{24m-3}$.

Now suppose that $S_x \cong {}^2F_4(q_0)$, where $q = q_0^r$ with $r$ prime. 
Writing $q_0 = 2^\ell$, we have
\[
|\mathcal{P}| = |S:S_x| = 2^{12\ell(r-1)} \frac{(2^{6r\ell}+1)(2^{4r\ell}-1)(2^{3r\ell}+1)(2^{r\ell}-1)}{(2^{6\ell}+1)(2^{4\ell}-1)(2^{3\ell}+1)(2^\ell-1)} < 2^{26\ell(r-1)+ 4}.
\]
However, Lemma~\ref{2part}(i) with $a=12\ell(r-1)$ gives $|\mathcal P| > 2^{36\ell(r-1)}$, a contradiction (because $\ell \ge 1$).

Finally, suppose that $S_x$ is as in one of the cases~(viii)--(x). 
Then the highest power of $2$ dividing $|S_x|$ is $2^5$ (arising in case~(ix)), so $|\mathcal{P}|=|S:S_x|$ is divisible by $2^{12m-5}$, and Lemma~\ref{2part}(i) therefore gives $|\mathcal{P}| > 2^{36m-15}$. 
On the other hand, we certainly have $|\mathcal{P}| < |S| < 2^{30m}$, which is a contradiction (because $36m-15 \le 30m$ if and only if $m \leq 5/2$, but we have $m \ge 3$).

\end{document}